\documentclass[11pt]{amsart}

\newtheorem{teo}{Theorem}[section]
\newtheorem{prop}[teo]{Proposition}
\newtheorem{cor}[teo]{Corollary}
\newtheorem{lema}[teo]{Lemma}

\theoremstyle{definition}
\newtheorem{dfn}[teo]{Definition}

\theoremstyle{remark}
\newtheorem{rem}[teo]{Remark}

\numberwithin{equation}{section}

\newcommand{\dist}{\ensuremath{\mathrm{dist} }}
\newcommand{\grad}{\ensuremath{\mathrm{grad}\ }}

\newcommand{\tub}{\ensuremath{\mathrm{Tub} }}

\newcommand{\F}{\ensuremath{\mathcal{F}}}
\newcommand{\E}{\ensuremath{\mathcal{E}}}
\newcommand{\planes}{\ensuremath{\mathcal{H}}}

\newcommand{\singularF}{\ensuremath{\mathcal{X}_{\mathcal{F}}}}

\newcommand{\dank}{\textsf{Acknowledgments.\ }}

\newcommand{\Holsing}{\ensuremath{\mathrm{Holsing}}}

\newcommand{\eop}{\hfill $\Box$\medskip}

\usepackage{amssymb}
\usepackage{amsmath}

\newcommand{\R}{\mbox{$\mathbf R$}}
\newcommand{\C}{\mbox{$\mathbf C$}}
\newcommand{\inn}[2]{\mbox{$\mathcal{h} #1,#2 \mathcal{i}$}}

\begin{document}

\title[S.r.f.s., transnormal maps ans basic forms]{ Singular riemannian foliations with sections, transnormal maps
and basic forms}

\author{Marcos M. Alexandrino \and Claudio Gorodski} 
 
\address{Marcos Martins Alexandrino \\
Instituto de Matem\'{a}tica e Estat\'{\i}stica,
Universidade de S\~{a}o Paulo (USP),
Rua do Mat\~{a}o, 1010, 
S\~{a}o Paulo, SP 05508-090, Brazil}

\email{malex@ime.usp.br}

\email{marcosmalex@yahoo.de}

\address{Claudio Gorodski \\
Instituto de Matem\'{a}tica e Estat\'{\i}stica,
Universidade de S\~{a}o Paulo (USP),
Rua do Mat\~{a}o, 1010, 
S\~{a}o Paulo, SP 05508-090, Brazil}

\email{gorodski@ime.usp.br}


\date{30th july 2006}

\thanks{The authors have been partially supported by FAPESP
and CNPq.}


\begin{abstract}

A singular riemannian foliation $\F$ on a complete riemannian 
manifold $M$ is said to admit sections if 
each regular point of $M$ is contained in a 
complete totally geodesic immersed 
submanifold $\Sigma$ that meets every leaf of $\F$ orthogonally 
and whose 
dimension is the codimension of the regular leaves of $\F$.

We prove that the algebra of basic forms of $M$ relative
to $\F$ is isomorphic 
to the algebra of those differential forms on $\Sigma$ that are 
invariant under the generalized Weyl pseudogroup of $\Sigma$. This 
extends a result of Michor for polar actions.
It follows from this result that the algebra of basic function 
is finitely generated if the sections are compact.

%A consequence of this result is a generalization of a 
%weel known result of G.W. Schwarz\ft{Need to change!} 
%stating that the algebra of basic functions is finitely generated. 

We also prove that the leaves of $\F$ coincide with the 
level sets of a transnormal map (generalization of isoparametric map) 
if $M$ is simply connected, the sections are flat and the leaves of 
$\F$ are compact.
This result extends previous results due to 
Carter and West, Terng, and Heintze, Liu and Olmos.

\end{abstract}

%%%%%%%%%%%%%%%%%%%%%%%%
%%%%%%%%%%% INTRODUCAO

\maketitle

\section{Introduction}

The main results of this paper are the following
two theorems. The first one generalizes 
previous results of Carter and West~\cite{CarterWest2},
Terng~\cite{terng} and Heintze, Liu and Olmos~\cite{HOL}
for isoparametric submanifolds. 
It can also be viewed as a converse to the main result 
in~\cite{Alex1}, and as a global version of one of the results 
in~\cite{Alex2}. 

\begin{teo}\label{teo-s.r.f.s-transnormal}
Let $\mathcal F$ be a singular riemannian foliation with sections
on a complete simply connected riemannian manifold $M$.
Assume that the leaves of $\mathcal F$ are compact and 
$\mathcal F$ admits a flat section of dimension $n$. 
Then the leaves of $\mathcal F$ are given by 
the level sets of a transnormal map $F:M\to\R^n$.
\end{teo}

The second theorem generalizes a result of 
Michor for basic forms relative to polar actions~\cite{Michor1,Michor2}, 
and will be used to prove Theorem~\ref{teo-s.r.f.s-transnormal}. 

\begin{teo}\label{teo-basic-forms}
Let $\mathcal F$ be a singular riemannian
foliation with sections on a complete riemannian
manifold $M$, and let $\Sigma$ be a section 
of $\mathcal F$. Then the immersion of
$\Sigma$ into $M$ induces an isomorphism 
between the algebra of basic differential
forms on $M$ relative to $\mathcal F$
and the algebra of differential forms 
on $\Sigma$ which are invariant 
under the generalized Weyl pseudogroup of $\Sigma$.
\end{teo}

Theorem~\ref{teo-basic-forms} of course includes the important case
of basic functions on $M$. In~\cite{Schwarz}
G.~Schwarz proved 
that the algebra of basic functions relative to the orbits of a
smooth action of a compact group on a compact manifold~$M$ 
is finitely generated. Using this result, Theorem~\ref{teo-basic-forms}
and a result of T\"{o}ben~\cite{Toeben}, we get the following consequence:

\begin{cor}\label{cor}
Let $\mathcal F$ be a singular riemannian
foliation with sections on a complete riemannian
manifold $M$. Assume that the sections of $\mathcal F$ are
compact submanifolds of $M$. Then the algebra
of basic functions on $M$ relative to
$\mathcal F$ is finitely generated.
\end{cor}

\emph{Singular riemannian foliations with sections} (\emph{s.r.f.s.},
for short) are singular riemannian foliations
in the sense of Molino~\cite{Molino} which admit 
transversal complete immersed manifolds that meet
all the leaves and meets them always orthogonally,
see section~2 for the definition.  
These were introduced by Boualem~\cite{Boualem}, 
and then by the first author~\cite{Alex1,Alex2}
as a simultaneous generalization of orbital foliations 
of polar actions of Lie groups (see e.g.~Palais and Terng~\cite{PTlivro}), 
isoparametric foliations in simply-connected space forms 
(see e.g.~Terng~\cite{terng}),
and foliations by parallel submanifolds of an equifocal submanifold 
with flat sections in a simply connected compact symmetric space 
(see e.g.~Terng and Thorbergsson~\cite{TTh1}).   
S.r.f.s.\ were further studied in~\cite{Alex3,Alex4,AlexToeben}, 
by T\"{o}ben~\cite{Toeben}, and also by 
Lytchak and Thorbergsson~\cite{LytchakThorbergsson}. 
By using suspensions of homomorphisms,
one can construct examples of s.r.f.s.\ with nonembedded or
exceptional leaves, and also inhomogeneous examples~\cite{Alex2}. 
Other techniques that are used to construct examples of 
s.r.f.s.\ on nonsymmetric spaces are surgery and
suitable changes of metric~\cite{AlexToeben}.  

An isoparametric submanifold in a simply-connected space form
can always be described as a regular level set of an isoparametric  
polynomial map~\cite{terng}. More generally, as proved by Heintze, Liu and Olmos~\cite{HOL}, 
an equifocal submanifold with flat 
sections in a simply connected compact symmetric space
can always be described as a regular level set 
of an analytic transnormal map
(a smooth map is called \emph{transnormal} if it is
an integrable riemannian submersion in a neighborhood
of any regular level set, see Definition \ref{dfn-transnormal-map}). 
In the case of s.r.f.s., recently 
it has been proved a local version of this
result, in that the plaques 
of a s.r.f.s.\ can always be described as level sets
of a locally defined transnormal map~\cite{Alex2};
Theorem~\ref{teo-s.r.f.s-transnormal} thus appears as the 
corresponding global statement. 
It is also worth noting that there is a global converse to the 
quoted result from~\cite{HOL}, namely,  
the regular leaves 
of an analytic transnormal map on a complete analytic
riemannian manifold are equifocal manifolds and leaves of a 
s.r.f.s.~\cite{Alex1}.

\medskip

Let $\mathcal F$ be a s.r.f.s.\ on a complete
riemannian manifold $M$. A smooth function on $M$
is called a \emph{basic function} if it is constant along the
leaves of $\mathcal F$. 
The definition of basic function 
can be extended to differential forms, in that
a differential form $\omega$ on $M$ is called a \emph{basic form}
if both $\omega$ and $d\omega$ vanish whenever at least one
of the arguments of $\omega$ (resp.~$d\omega$)
is a vector tangent to a leaf of $\mathcal F$.
Palais and Terng~\cite{PT} considered the case of a 
polar action of a Lie group $G$ on a complete riemannian manifold $M$
and proved that the restriction from $M$ to a section $\Sigma$ 
induces an isomorphism between the algebra of basic 
functions on $M$ relative to the orbital foliation
and the algebra of functions on $\Sigma$ which are invariant
under the generalized Weyl group of $\Sigma$. 
Michor~\cite{Michor1,Michor2} extended Palais and Terng's result
to basic forms relative to a polar action. 
In this context,
Theorem~\ref{teo-basic-forms} is a generalization of these results to 
s.r.f.s.\/. Here it is important to remark that
a s.r.f.s.\ admits a generalized Weyl pseudogroup 
which acts on a section, but, in general, this is not a 
Weyl group, see section~2. 

We finish this introduction with some remarks about  
Wolak's claim to have proven Theorem~\ref{teo-basic-forms} in~\cite{Wolak}
under the additional hypothesis that the leaves
be compact. In our opinion, there are two  problems with 
his arguments. The first one is that he has used a Weyl
pseudogroup but has not defined it properly. It would appear
that he has used the pseudogroup constructed by Boualem
in~\cite{Boualem}. 
Even if this is the case, Boualem's pseudogroup is 
often smaller than the needed pseudogroup, which is 
correctly defined in~\cite[Definition~2.6]{Alex2} (see also 
remarks in~\cite{Alex2,AlexToeben}). In fact,
in order to properly define the Weyl pseudogroup, 
one needs the equifocal property or something equivalent to it. 
The second problem that we found with Wolak's arguments 
is related to a property of s.r.f.s.. 
He incorrectly claimed at the end of Proposition~1 in~\cite{Wolak} that the 
restriction of the foliation to a slice must be homogeneous. 
This claim is false since there exist many examples of 
isoparametric foliations with inhomogenous leaves in 
euclidean space~\cite{FKM}.

\medskip

\dank The authors are grateful to Professor Gudlaugur Thorbergsson  
for useful suggestions.

%%%%%%%%%%%%%%%%%%%%%%%%%%%%%%%%%%%%%%%%%%%%%%%%%%%%%%%%%%%%%%%%%%%%%%%%%%%%%%%%%%%%%%%%%%%%%%%%
%%%%%%%%%FACTS ABOUT SRFS
%%%%%%%%%%%%%%%%%%%%%%%%%%%%%%%%%%%%%%%%%%%%%%%%%%%%%%%%     

\section{Facts about s.r.f.s.}

In this section, we recall some results about s.r.f.s.\ that will be used 
in this text. Details can be found in~\cite{Alex2,AlexToeben}. 
Throughout this section, we assume that $\F$ is a 
singular riemannian foliation with sections
on a complete riemannian manifold $M$; 
we start by recalling its definition.

\begin{dfn}
A partition $\F$ of a complete riemannian manifold $M$ by connected 
immersed submanifolds (the \emph{leaves}) is called a \emph{singular 
riemannian foliation with sections} of $M$ (\emph{s.r.f.s.}, for short) if it 
satisfies
the following conditions:
\begin{enumerate}
\item $\F$ is \emph{singular},
i.e.~the module $\singularF$ of smooth vector fields on $M$ that are 
tangent at each point to the corresponding leaf acts transitively on each 
leaf. In other words, for each leaf $L$ and each $v\in TL$ with 
footpoint $p,$ there exists $X\in \singularF$ with $X(p)=v$.
\item  The partition is \emph{transnormal}, i.e.~every 
geodesic that is perpendicular to a leaf at one point remains 
perpendicular to every leaf it meets.
\item For each regular point $p$, the set $\Sigma :=\exp_{p}(\nu_p 
L_{p})$ 
is a complete immersed submanifold that meets all the leaves and meets
them always orthogonally. The set  
$\Sigma$ is called a \emph{section}.
\end{enumerate}
\end{dfn}

\begin{rem}
In \cite{Boualem} Boualem dealt with a 
singular riemannian foliation $\F$ on a complete manifold $M$ 
such that the distribution of normal spaces of the regular leaves  
is integrable. It was proved in \cite{Alex4} that such an 
$\F$ must be a s.r.f.s.\/ and, in addition,
the set of regular points is open and dense in each section. 
\end{rem}

A typical example of a s.r.f.s is the partition formed by parallel submanifolds  of an isoparametric submanifold $N$ of an euclidean space.
A submanifold $N$ of an euclidean space is called \emph{isoparametric} 
if its normal bundle is flat and the principal curvatures along any 
parallel 
normal vector field are constant. 
Theorem~\ref{sliceteorema} below shows how s.r.f.s.\ and 
isoparametric foliations are related to each other.
In order to state this theorem, we need the concepts of slice and  
local section. 
Let $q\in M$, and let $\tub(P_{q})$ be a tubular neighborhood of a 
plaque $P_{q}$ that contains $q$. Then the connected component of 
$\exp_{q}(\nu P_{q})\cap \tub(P_{q})$  that contains $q$ is called a 
\emph{slice} at $q$ and is usually denoted by $S_{p}.$ 
A \emph{local section} $\sigma$ (centered at $q$) of a section $\Sigma$  is a connected component $\tub(P_{q})\cap\Sigma$. 
 
\begin{teo}[\cite{Alex2}]
\label{sliceteorema}
Let $\F$ be a s.r.f.s. on a complete riemannian manifold $M.$ 
Let $q$ be a singular point of $M$ and let 
$S_{q}$ a slice at $q.$ Then
\begin{enumerate}
\item Let $\epsilon$ be the radius of the slice $S_{q}.$ Denote  
$\Lambda(q)$  the set of 
local sections $\sigma$ containing $q$ such that $\dist(p,q)<\epsilon$ 
for each $p\in\sigma.$ 
 Then  $S_{q}= \cup_{\sigma\in\Lambda (q)}\, \sigma$.
\item $S_{x}\subset S_{q}$ for all $x\in S_{q}$.
\item $\F|S_q$ is a s.r.f.s.\ on $S_{q}$ with the induced metric from 
$M$.
\item $\F|S_q$ is diffeomorphic to an isoparametric foliation on an 
open 
subset of $\mathbf{R}^{n}$, where $n$ is the dimension of $S_{q}$.
\end{enumerate}
\end{teo}

From (d), it is not difficult to derive the following corollary.

\begin{cor}
\label{estratificacao-singular}
Let $\sigma$ be a local section. Then the set of singular points of 
$\F$ 
that are contained in $\sigma$ is a finite union of totally 
geodesic hypersurfaces. 
These hypersurfaces are mapped by a diffeomorphism to the focal 
hyperplanes 
contained in a section of an isoparametric foliation on an open subset 
of an euclidean space.   
\end{cor}

We will call the set of singular points of $\F$ contained in $\sigma$ 
the \emph{singular stratification of the local section} $\sigma$. 
Let $M_{r}$ denote the set of regular points in $M.$ A \emph{Weyl 
Chamber} of 
a local section $\sigma$ is the closure in $\sigma$ of a connected 
component 
of $M_{r}\cap\sigma$. One can prove that a Weyl Chamber of a 
local section is a convex set.  

It also follows from Theorem~\ref{sliceteorema} that 
\emph{the plaques of a s.r.f.s.\ are always level sets of a transnormal map}, 
whose definition we recall now. 

\begin{dfn}[Transnormal  Map]
\label{dfn-transnormal-map}
Let  $M^{n+q}$  be a complete riemannian manifold. A smooth map   
$F=(f_{1},\ldots, f_{q}):M^{n+q}\rightarrow \mathbf{R}^{q}$ is 
called a \emph{transnormal map} if the following assertions hold:
\begin{enumerate}
\item[(0)]  $F$ has a regular value.
\item[(1)]  For every regular value $c$, there exists a neighborhood $V$ 
of $F^{-1}(c)$ in $M$ and smooth functions $b_{i\,j}$ on $F(V)$ such that, for every $x\in V,$ 
$\langle\grad f_{i}(x),\grad f_{j}(x)\rangle= b_{i\,j}\circ F(x).$
\item[(2)]  There is a sufficiently small neighborhood of each regular 
level set on which, for every $i$ and $j$, the bracket
$[\grad f_{i},\grad f_{j}]$ is a linear combination of 
$\grad f_{1},\ldots,\grad f_{q}$, where the coefficients are 
functions of $F$.  
\end{enumerate}
\end{dfn}
This definition is equivalent to saying that the map $F$ has a regular value and for each regular value $c$ there exists a neighborhood  $V$ of $ F^{-1}(c)$ in $M$ such that $F \mid_{V}\rightarrow F(V)$ is an integrable riemannian submersion, where the metric $(g_{i\,j})$ of  $F(V)$ is the inverse matrix of $(b_{i\,j}).$
In particular, a transnormal map $F$ is said to be an \emph{isoparametric map} if  $V$ can be chosen to be $M$  and   $\bigtriangleup f_{i}= a_{i}\circ F,$ where $a_{i}$ are smooth functions. 
As we have remarked in the indroduction, each isoparametric submanifold in an euclidian space
can always be described as a regular level set of an isoparametric  
polynomial map (see \cite{terng} or \cite{PTlivro}).

In \cite{TTh1}, Terng and Thorbergsson introduced the concept of 
equifocal submanifolds with flat sections in symmetric spaces in 
order to generalize the definition of isoparametric submanifolds 
in euclidean space. 
Next we review the slightly more general
definition of equifocal submanifolds in riemannian 
manifolds.

\begin{dfn}
\label{dfn-equifocal}
A connected immersed submanifold $L$ of a complete riemannian manifold 
$M$ is called \emph{equifocal} if it satisfies the following 
conditions:
\begin{enumerate}
\item The normal bundle $\nu(L)$ is flat.
\item $L$ has sections, i.e.~for each~$p\in L$, 
the set $\Sigma :=\exp_{p}(\nu_p L_{p})$ is a complete immersed totally 
geodesic submanifold.
\item For each parallel normal field $\xi$ on a neighborhood $U \subset L$,  
the derivative of the map $\eta_{\xi}:U\to M$ defined 
by $\eta_{\xi}(x):=\exp_{x}(\xi)$ has constant rank.
\end{enumerate}
\end{dfn}

The next theorem relates s.r.f.s.\ and equifocal submanifolds.

\begin{teo}[\cite{Alex2}]
\label{frss-eh-equifocal}
Let $L$ be a regular leaf of a s.r.f.s.\ $\F$ of a 
complete riemannian manifold $M$.
\begin{enumerate}
\item Then $L$ is equifocal. In particular, the union of the 
regular leaves that 
have trivial normal holonomy is an open and dense set in $M$ provided 
that all the leaves are compact.
\item Let $\beta$ be a smooth curve of $L$ and $\xi$ a parallel normal 
field to
$L$ along~$\beta$. 
Then the curve $\eta_{\xi}\circ \beta$ belongs to a leaf of $\F$.
\item Suppose that $L$ has trivial holonomy and let $\Xi$ denote the 
set of 
all parallel normal fields on $L$. 
Then $\F=\{\eta_{\xi}(L)\}_{\xi\in \, \Xi}.$ 
\end{enumerate}
\end{teo}

The above theorem allows us to define the 
singular holonomy map, which will be very useful to study $\F$.  

\begin{prop}[Singular holonomy]
\label{prop-holonomia-singular}
Let $L_{p}$ be a regular leaf, let $\beta$ be a smooth curve in~$L_{p}$ 
and let $[\beta]$ denote its homotopy class. 
Let~$U$ be a local section centered at $p=\beta(0)$. 
Then there exists a local section $V$ centered at $\beta(1)$ and 
an isometry
$\varphi_{[\beta]}:U\to V$ with the following properties:
\begin{enumerate}
\item[1)]$\varphi_{[\beta]}(x)\in L_{x}$ for each $x\in U$.
\item[2)]$d\varphi_{[\beta]}\xi(0)=\xi(1),$ where $\xi$ is a  parallel 
normal field along $\beta$.
\end{enumerate}
\end{prop}

An isometry as in the above proposition
is called the \emph{singular holonomy map along $\beta$}. 
 
We remark that, in the definition of the singular holonomy map, 
singular points can be contained in the domain $U.$  
If the domain $U$ and the range $V$ are sufficiently small, then the 
singular holonomy map coincides with the usual holonomy map along 
$\beta$.

Theorem \ref{sliceteorema} establishes a relation between s.r.f.s.\ and 
isoparametric foliations. Similarly as in 
the usual theory of isoparametric submanifolds, it is 
natural to ask if we can define a (generalized) Weyl group action on 
$\sigma$. The following definitions and results 
deal with this question.

\begin{dfn}[Weyl pseudogroup $W$]
\label{definitionWeylPseudogroup}
 The pseudosubgroup generated by all singular holonomy maps 
$\varphi_{[\beta]}$ such that $\beta(0)$ and $\beta(1)$ belong to the 
same 
local section~$\sigma$ is called the \emph{generalized Weyl 
pseudogroup}
of~$\sigma$. Let $W_{\sigma}$ denote this pseudogroup. 
In a similar way, we define $W_{\Sigma}$ for a section $\Sigma$. 
Given a slice~$S$, we define~$W_{S}$ as the set of all singular 
holonomy 
maps~$\varphi_{[\beta]}$ such that~$\beta$ is contained in the 
slice~$S$.
\end{dfn}

\begin{rem}
Regarding the definition of pseudogroups and orbifolds,
see Salem~\cite[Appendix D]{Molino}. 
\end{rem}

\begin{prop}
\label{propWisinvariant}
Let $\sigma$ be a local section. Then the reflections in the 
hypersurfaces of the singular stratification of the local section $\sigma$ leave 
$\F|\sigma$  invariant. Moreover these reflections are elements of 
$W_{\sigma}.$
\end{prop}

By using the technique of suspension, 
one can construct an example of a s.r.f.s.\ 
such that  $W_{\sigma}$ is larger than the pseudogroup 
generated by the reflections in the hypersurfaces of the 
singular stratification of $\sigma$. On the other hand,
a sufficient condition 
to ensure that both pseudogroups coincide is that the leaves of $\F$ 
have trivial normal holonomy and be compact.
So it is natural to ask under which conditions we can garantee that the normal holonomy of regular leaves are trivial. The next result is 
concerned with this question.  

\begin{teo}[\cite{AlexToeben}]
\label{teo-fundamental-domain}
Let $\F$ be a s.r.f.s.\ on a simply connected riemannian manifold $M$. 
Suppose also that the leaves of $\F$ are compact. Then
\begin{enumerate}
\item Each regular leaf has trivial holonomy.
\item $M/ \F$ is a simply connected Coxeter orbifold.
\item Let $\Sigma$ be a section of $\F$ and 
let $\Pi:M\rightarrow M/\F$ be the canonical projection. 
Denote by $\Omega$ a connected component of the set of regular points 
in 
$\Sigma$. Then  $\Pi:\Omega\rightarrow M_{r}/\F$ and 
$\Pi:{\overline{\Omega}}\to M/\F$ are homeomorphisms, where $M_r$ 
denotes the set of regular points in $M$. In addition, 
$\Omega$ is convex, i.e.~for any two points~$p$ and~$q$ in~$\Omega$, 
every minimal geodesic segment between $p$ and $q$ lies 
entirely in $\Omega$.
\end{enumerate}
\end{teo}

%%%%%%%%%%%%
%%%%%%%%%%%%%%%%%%%%%%%%%%%
%PROVA DO TEO FORMAS BASICAS
%%%%%%%%%%%%%%%%%%%%%%%%%%%%%%%%%%%%%%%%%%%%%%%%%%%%%%%%%%%%%%%%%%%%%%%%
%%%%%%%%%%%
\section{Proof of Theorem \ref{teo-basic-forms}}

Throughout this section we assume that $\F$ is a s.r.f.s.\ 
on a complete riemannian manifold $M$ and 
prove Theorem~\ref{teo-basic-forms}. 
We start by recalling the definition of basic forms.

\begin{dfn}[Basic forms]
A differential $k$-form $\omega$ is said to be \emph{basic} if, 
for all $X\in \singularF$, we have:
\begin{enumerate}
\item[(a)] $i_{X}\omega=0,$
\item[(b)] $i_{X}d\omega=0.$
\end{enumerate}
\end{dfn}

In the course of the proof of the theorem, we will also need the 
concept of  differential form invariant by holonomy. 
As soon as Theorem \ref{teo-basic-forms} is proved, 
it will be clear that these two concepts are in fact equivalent.

\begin{dfn}\label{holon-inv-form}
A differential $k$-form $\omega$ is said to be
\emph{invariant by holonomy} 
if:
\begin{enumerate}
\item[(a)] $i_{X}\omega=0$ for all $X\in \singularF.$ 
\item[(b)] Let $\sigma$ and $\tilde{\sigma}$ be local sections 
and $\varphi:\sigma\rightarrow \tilde{\sigma}$ a singular holonomy map. 
Let $I:\sigma\rightarrow M$ and $\tilde{I}:\tilde{\sigma}\rightarrow M$ 
be
the inclusions of $\sigma$ and $\tilde{\sigma}$ in $M$. Then
$\varphi^{*}(\tilde{I}^{*}\omega)=I^{*}\omega$
\end{enumerate}
A $k$-form $\omega$ is said to be 
\emph{invariant by regular holonomy} if it satisfies 
the definition above with the condition that 
the map $\varphi:\sigma\to\tilde{\sigma}$
can be only a regular holonomy. 
\end{dfn}

\begin{lema}
\label{lema-holinv-eh-basica}
\begin{enumerate}
\item[(a)] Forms invariant by holonomy are basic forms.
\item[(b)] Basic forms are  invariant by regular holonomy.
\end{enumerate}
\end{lema}
\begin{proof}
(a) Let $\omega$ be a form invariant by holonomy and let 
$X\in \singularF.$  We want to prove that
\begin{equation}
\label{lema-holinv-eh-basica-Eq0}
i_{X}d\omega=0.
\end{equation}

First we prove this equation 
for regular points. Let $p$ be a regular point and let
$P_p$ be a plaque of $\F$ that contains $p$.
Using Theorem \ref{frss-eh-equifocal} and the normal exponential map 
$\exp^{\nu}: \nu(P_p)\rightarrow M$, we can construct a vector 
field $\tilde{X}$ on a neigborhood of $p$ such that
 \begin{enumerate}
 \item $\tilde{X}(y)=X(y)$ for $y\in P_{p}.$
 \item $\varphi_{t}:=\psi_{t}|_{\sigma_0}$ is a regular
holonomy map, where
$\psi_t$ is the flow of $\tilde{X}$ and $\sigma_0$ 
is a local section that contains $p$.
\end{enumerate}

Define $\sigma_t:=\varphi_t(\sigma_0)$, and let 
 $I_t:\sigma_t\rightarrow M$ denote  the inclusion of $\sigma_t$ in 
$M$.
 Since $\omega$ is  invariant by holonomy, we have
\begin{equation}
\label{lema-holinv-eh-basica-Eq1}
I^{*}_{0}\psi_{t}^{*}\omega=\varphi_{t}^{*}I_{t}^{*}\omega = 
I_{0}^{*}\omega.
\end{equation}
On the other hand, it follows from $\varphi_{t}(x)\in L_x$ that
\begin{equation}
\label{lema-holinv-eh-basica-Eq2}
i_{Y(p)}\psi^{*}_{t}\omega=0
\end{equation}
for each $Y\in\singularF$.
Putting together
Equations~(\ref{lema-holinv-eh-basica-Eq1}) 
and~(\ref{lema-holinv-eh-basica-Eq2}), we get that 
\begin{equation}
\label{lema-holinv-eh-basica-Eq3}
\psi_{t}^{*}\omega=\omega
\end{equation}
for small $t.$

Equation (\ref{lema-holinv-eh-basica-Eq3}) and the definition of $\tilde{X}$ 
now yield that 
\begin{eqnarray*}
0&=& L_{\tilde{X}}\omega(p)\\
 &=& i_{\tilde{X}(p)}d\omega+ 
d(i_{\tilde{X}}\omega)(p)\\
 &=&i_{X(p)}d\omega.
\end{eqnarray*} 
We have shown that Equation (\ref{lema-holinv-eh-basica-Eq0}) holds on the 
set of regular points of $\F$. 
Since this set is dense in $M$,
this finishes the proof. 

(b) Let $\omega$ be a basic form and let 
$\varphi_{[\beta]}:\sigma_0\rightarrow\sigma_1$ 
be a regular holonomy, where $\sigma_0$ and $\sigma_1$ 
are local sections that contains only regular points and $\beta$ 
is a curve contained in a regular leaf such that $\beta(0)\in\sigma_0$ 
and  $\beta(1)\in\sigma_1.$ Let $0=t_0<\cdots<t_n=1$ be a partition 
such that $\beta_i:=\beta|_{[t_{i-1},t_{i}]}$ is a curve contained 
in a distinguished neighborhood $U_i$, i.e.~the plaques of $\F$ in 
$U_i$ are fibers of a submersion. 
Since $\varphi_{\beta}=\varphi_{\beta_n}\circ\ldots\circ\varphi_{\beta_1}$, 
in order to see that $\omega$ is invariant by regular holonomy, 
it suffices to prove that 
\begin{eqnarray} 
\label{Eq-lema-basic-eh-regular-invariant-form}
\varphi_{[\beta_i]}^{*}(I_{i}^{*}\omega)=I_{i-1}^{*}\omega, 
\end{eqnarray}
where $I_{i}:\sigma_i\rightarrow M$ is the inclusion in $M$ 
of a local section $\sigma_i$ that contains $\beta(t_i).$

Since $\beta_i$ is contained in a distinguished neighborhood, 
we can construct a field $X$ such that 
$\varphi_{[\beta_i]}=\psi_{t_{i}}|_{\sigma_{i-1}},$ where 
$\psi$ is the flow of $X$. 
The fact the $\omega$ is basic gives that $L_{X}\omega=0$. 
Therefore $\psi_{t}^{*}\omega=\omega$, and 
this implies Equation~(\ref{Eq-lema-basic-eh-regular-invariant-form}).
\end{proof}

Let $I:\Sigma\rightarrow M$ be the immersion of the section $\Sigma$ in 
$M$. We divide the proof of Theorem~\ref{teo-basic-forms} 
into the following two lemmas.

\begin{lema}
\label{lema-I-eh-injetora}
The map $I^{*}: \Omega^{k}_{b}(M)\rightarrow 
\Omega^{k}(\Sigma)^{W_\Sigma}$ is well defined and injective, 
where $\Omega^{k}_{b}(M)$ denotes the algebra of basic $k$-forms 
on $M$ and $\Omega^{k}(\Sigma)^{W_\Sigma}$ denotes the algebra of  
$W_\Sigma$-invariant $k$-forms on $\Sigma$.
\end{lema}
\begin{proof}
Let $\omega\in \Omega^{k}_{b}(M)$.
We first verify that 
\begin{eqnarray}
\label{Eq-1-injetividade-da-inclusao}
 I^{*}\omega \in \Omega^{k}(\Sigma)^{W_\Sigma}.
\end{eqnarray}
  In fact, Lemma~\ref{lema-holinv-eh-basica}(b) implies that
\begin{eqnarray}
\label{Eq-2-injetividade-da-inclusao}
\varphi^{*}(I^{*}\omega)|_{\Sigma\cap M_{r}}=I^{*}\omega|_{\Sigma\cap M_{r}},
\end{eqnarray}
for a singular holonomy $\varphi\in W_{\Sigma}.$
Equation (\ref{Eq-1-injetividade-da-inclusao}) then
follows from~(\ref{Eq-2-injetividade-da-inclusao}), 
as the set $M_r\cap\Sigma$ is dense in $\Sigma$. 
To show that $I^{*}$ is injective, we need to see
that $I^{*}\omega=0$ implies $\omega=0$. 
This follows again from~Lemma \ref{lema-holinv-eh-basica}(b) 
and the denseness of $M_r$ in $M.$ 
\end{proof}

\begin{lema}
\label{lema-I-eh-sobrejetora}
The map $I^{*}:\Omega^{k}_{b}(M)\to\Omega^{k}(\Sigma)^{W_\Sigma}$ 
is surjetive.
\end{lema}

\begin{proof}
According to Lemma~\ref{lema-holinv-eh-basica}(a), 
it suffices to check that
\[ I^{*}: \Omega^{k}(M)^{\Holsing}\to
\Omega^{k}(\Sigma)^{W_\Sigma}\]
is surjective, where 
$\Omega^{k}(M)^{\Holsing}$ is the algebra of differential $k$-forms on $M$
that are invariant by holonomy.

Let $\omega\in \Omega^{k}(\Sigma)^{W_\Sigma}$.
We will construct $\tilde{\omega}\in \Omega^{k}(M)^{\Holsing}$ such 
that $I^{*}\tilde{\omega}=\omega$.
Let $\tilde{q}$ be a point of $M$. 
To begin with, we set $i_{Y}\tilde{\omega}=0$ for all 
$Y\in T_{\tilde{q}}P_{\tilde{q}}$, where $P_{\tilde{q}}$ 
is the plaque that contains $\tilde{q}$.
Next, we must define $\tilde{\omega}|_{T_{\tilde{q}}S_{\tilde{q}}}$,
where $S_{\tilde{q}}$ is a slice at $\tilde{q}$. 

Suppose first that $\tilde{q}$ is a regular point. 
In this case, the slice $S_{\tilde{q}}$ 
is a local section. 
Let $\sigma\subset\Sigma$ be a local section that contains 
a point $q\in L_{\tilde{q}}\cap\Sigma$ and let $\varphi:\sigma\rightarrow 
S_{\tilde{q}}$ be a singular holonomy.
We define $\tilde{\omega}|_{S_{\tilde{q}}}=(\varphi^{-1})^{*}\omega.$ 
This definition does not depend on $\sigma$ and $\varphi$ since
$\omega\in \Omega^{k}(\Sigma)^{W_\Sigma}.$

Next, suppose that $\tilde{q}$ is a singular 
point. In this case, the slice $S_{\tilde{q}}$ is no longer a local section, 
but, on the contrary, it is the union of the local sections that contain 
$\tilde{q}$ (see Theorem~\ref{sliceteorema}(a)). 
This leads us to consider the intersection of all those local sections 
that contain $\tilde{q}$.
Denote by $T$ be the connected component of the minimal  
stratum of the foliation $\F\cap S_{\tilde{q}}$ that 
contains $\tilde{q}$. It follows from Theorem~\ref{sliceteorema} 
and from the theory of isoparametric submanifolds 
that $T$ is a union of singular points of the foliation 
$S_{\tilde{q}}\cap\F$, and $T$ is the intersection of 
the local sections that contain $\tilde{q}$. 

%Note that if the point $\tilde{q}$ belonged to the section $\Sigma$ then the fact that $\omega\in \Omega^{k}(\Sigma)^{W_\Sigma}$ would imply that 
%$i_{Y} \omega=0$ if $Y(\tilde{q})\in \nu(T).$ 

In order to motivate the next step in the 
construction of $\tilde\omega$, we remark that 
for a point $\tilde q\in\Sigma$ and $Y(\tilde{q})\in \nu_{\tilde{q}}(T)$,
we have that $i_{Y} \omega=0$. For the purpose of checking
this remark, consider a local section $\sigma$ 
such that $Y\in T_{\tilde{q}}\sigma$. 
Then the theory of isoparametric submanifolds and Theorem~\ref{sliceteorema} 
allow us to choose a basis $\{Y_i\}$ of 
$\nu_{\tilde{q}} T \cap T_{\tilde{q}}\sigma$  
such that each $Y_i$ is orthogonal to a wall $H_i$. 
Due to the invariance of $\omega$ under the reflection on $H_i$, 
we have that $i_{Y_i}\omega|_{H_i}=i_{-Y_i}\omega|_{H_i}$.
Hence $i_{Y_i} \omega=0$. 
Since $\{Y_i\}$ is a basis, we conclude that $i_{Y} \omega=0$, as 
wished.

%To check the above claim we note  that $i_{Y} \omega=0$ for a $Y$ orthogonal to a  wall $H$. 
%In fact let $\varphi$ be a reflection in $H.$ Then the fact that   $\omega\in \Omega^{k}(\Sigma)^{W_\Sigma}$ implies that %$i_{Y}\omega|_{H}=i_{-Y}\omega|_{H}$ and hence $i_{Y} \omega=0.$ Now the theory of isoparametric submanifolds and  Theorem \ref{sliceteorema} imply %that we can choose a basis $\{Y_i\}$ of $\nu_{\tilde{q}} T$  such that each $Y_i$ is  orthogonal to a wall.

Based on the previous remark, we set
$i_{Y}\tilde{\omega}=0$ for $Y(\tilde{q})\in \nu(T)$ 
and arbitrary $\tilde q\in M$. 

It remains to define $\tilde{\omega}|_{T}$. To do that,
choose a local section 
$\tilde{\sigma}$ that contains $\tilde{q}$, 
a point $q\in L_{\tilde{q}}\cap\Sigma$, 
a local section $\sigma\subset\Sigma$ that contains $q$, and 
a holonomy $\varphi:\sigma\rightarrow\tilde{\sigma}$. 
Then we set
$\tilde{\omega}|_{T}=(\varphi^{-1})^{*}\omega|_{T}$. 
Let us show that 
the definition does not depend on $\sigma$, $\tilde{\sigma}$ and 
$\varphi$ by using that $\omega\in \Omega^{k}(\Sigma)^{W_\Sigma}$.
Indeed, for $i=1,2$, let $\sigma_i$ be a local section that 
contains a point $q_{i}\in L_{\tilde{q}}\cap\Sigma$, let  
$\tilde{\sigma}_i$ be a local section of $\tilde{q}$,
and let $\varphi_i:\sigma_{i}\rightarrow\tilde{\sigma}_{i}$ be
a singular holonomy. 
Denote by 
$\varphi_{2\,1}:\tilde{\sigma}_{1}\to\tilde{\sigma}_{2}$ 
a singular holonomy in $W_{S_{\tilde{q}}}$ (see Definition \ref{definitionWeylPseudogroup}). Then 
Theorem~\ref{sliceteorema} and the theory of 
isoparametric submanifolds force $\varphi_{2\,1}|_{T}$ 
to be the identity. Owing to our assumption on $\omega$,
$(\varphi^{-1}_{2}\circ\varphi_{2\,1}\circ \varphi_{1})^{*}\omega=\omega$.
Hence $(\varphi^{-1}_{1})^{*}\omega|_{T}=(\varphi^{-1}_{2})^{*}\omega|_{T}$.

We have already constructed the $k$-form $\tilde{\omega}$. 
It also follows from the construction that $\tilde{\omega}$ 
is invariant by holonomy in the sense that it satisfies the 
conditions~(a) and~(b) in Definition~\ref{holon-inv-form}. 
Now it only remains to prove that $\tilde{\omega}$ is 
smooth. It suffices to prove that 
in a neighborhood of an arbitrary point~$\tilde{q}$.  

If $\tilde{q}$ is a regular point, then there exists only one (germ of) 
local section $\sigma$ that contains $\tilde{q}$.
By construction, $\tilde{\omega}|_{\sigma}$ is smooth. 
Since $\tilde{\omega}$ is invariant by holonomy, 
we deduce that $\tilde{\omega}$ is smooth in a distinguished 
neighborhood of~$\tilde{q}$.

Next, we suppose that $\tilde{q}$ is a singular point. 
Let $\psi:S_{\tilde{q}}\to U\subset\textbf{R}^{n}$ be the 
diffeomorphism that sends the s.r.f.s.~$\F\cap S_{\tilde{q}}$ 
into an isoparametric foliation~$\hat{\F}$ in the open set~$U$ 
of euclidean space~$\textbf{R}^{n}$ 
(where~$n$ is the dimension of~$S_{\tilde{q}}$).
Note that $\tilde{\omega}|_{S_{\tilde{q}}}$ is invariant by each 
holonomy $\varphi_{[\beta]}$, where $\beta$ is a curve 
contained in the slice $S_{\tilde{q}}$. 
Since the diffeomorphism $\psi$ sends local sections into local sections, 
we conclude that 
$(\psi^{-1})^{*}( \tilde{\omega}|_{S_{\tilde{q}}})$ is  
invariant by the holonomy of the foliation $\hat{\F}$. 

Fix a section $V$ of the isoparametric foliation 
$\hat{\F}$, set $\dim V=l$, and select a
minimal set of homogeneous generators 
$\kappa_1,\cdots, \kappa_{l}$ of
the algebra $\mathbf{R}[V]^{\hat{W}}$ of $\hat{W}$-invariant functions 
of $V$, where $\hat{W}$ is the Coxeter group of the isoparametric 
foliation $\hat{\F}$ (see e.g.~\cite{PT}). 

By construction, we have that $\tilde{\omega}$ restricted to a 
local section that contains $\tilde{q}$ is smooth. 
Since $\psi^{-1}(V)$ is a local section, 
$\tilde{\omega}|_{\psi^{-1}(V)}$ is smooth. 
Therefore  $(\psi^{-1})^{*} (\tilde{\omega}|_{S_{\tilde{q}}})|_{V}$
is smooth and invariant by $\hat{W}$. 
Then it follows as in 
Michor~\cite[Lemma 3.3 and proof of Theorem 3.7]{Michor1} that
\[
(\psi^{-1})^{*} (\tilde{\omega}|_{S_{\tilde{q}}})|_{V}=
\sum \eta_{i_{1}\cdots 
i_{j}} d \kappa_{i_1}\wedge \cdots\wedge d \kappa_{i_j},
\]
where $\eta_{i_{1}\cdots i_{j}}$ is a smooth $\hat{W}$-invariant function.
In view of Schwarz~\cite{Schwarz}, we can write   
$\eta_{i_{1}\cdots i_{j}}= \lambda_{i_{1}\cdots i_{j}}(\kappa_1,\ldots,\kappa_l)$ for a smooth 
function $\lambda_{i_{1}\cdots i_{j}}$ on $\mathbf{R}^{l}$.
By~\cite{PTlivro}, we can extend each $\kappa_i$ to a
$\hat{\F}$-invariant function $\hat\kappa_i$ on~$U$.
%Since $(\psi^{-1})^{*} (\tilde{\omega}|_{S_{\tilde{q}}} )\in 
%\Omega^{k}(U)^{\hat{\F}}$ we have
Since $(\psi^{-1})^{*}( \tilde{\omega}|_{S_{\tilde{q}}})$ is  
invariant by the holonomy of the foliation $\hat{\F}$, we have 
\[
(\psi^{-1})^{*} (\tilde{\omega}|_{S_{\tilde{q}}})=\sum  \lambda_{i_{1}\cdots 
i_{j}}(\hat\kappa_1,\ldots,\hat\kappa_l)\,d \hat{\kappa}_{i_1} \wedge\cdots \wedge d \hat{\kappa}_{i_j}.
\]
Therefore
\[
\tilde{\omega}|_{S_{\tilde{q}}}=\sum  \lambda_{i_{1}\cdots 
i_{j}}(f_1,\ldots,f_l)\, d f_{i_1} \wedge \cdots \wedge d f_{i_j},
\]
where we have set $f_i= \hat{\kappa}_i\circ\psi$.
This equation already shows that 
$\tilde{\omega}|_{S_{\tilde{q}}}$ is smooth on 
$S_{\tilde{q}}$. Finally, extend each $f_i$ 
to a function defined on a tubular 
neigborhood of $\tilde{q}$ and denoted by the same
letter by setting it to be
constant along each plaque in that neigborhood.
Then the preceding equation and the invariance by holonomy imply that 
\[
\tilde{\omega}=\sum  \lambda_{i_{1}\cdots 
i_{j}}(f_1,\ldots,f_l)\, d f_{i_1} \wedge \cdots \wedge d f_{i_j},
\]
on a tubular neighborhood of $\tilde{q}$, which finally 
shows that 
$\tilde{\omega}$ is smooth in a neighborhood of $\tilde\omega$. 
Hence $\tilde\omega$ is smooth. 
\end{proof}

%%%%%%%%%%%%%%%%%%%%%%%%%%%%%%%%%%%%%%%%%%%%%%%%%%%%%%%%%%%%%%%%%%%%
%%%% COROLARIO
%%%%%%
%%%%%%%%%%%%%%%%%%%%%%%%%%%%%%%%%%%%%%%%%%%%%%%%%%%%%%%%%%%%%%%%%%%%%%%

\section{Proof of Corollary~\ref{cor}}

We rewrite the statement of Corollary~\ref{cor} as follows. 

\begin{cor}
\label{cor-basic-basis}
Let $\F$ be a s.r.f.s.\ on a complete riemannian manifold~$M$, and 
let $\E(M)^{\F}$ denote the space of basic functions.  
Suppose that the sections of $\F$ 
are compact submanifolds of $M$. Then there exist functions 
$f_1,\ldots,f_n\in \E (M)^{\F}$ such that
\[F^{*}\E(\mathbf{R}^{n})=\E (M)^{\F},\]
where $F=(f_1,\ldots,f_n)$.
\end{cor}
\begin{proof}
According to T\"{o}ben~\cite{Toeben}, 
there exists a section~$\Sigma$ 
such that the Weyl pseudogroup~$W_{\Sigma}$ is in fact a group. 
By Ascoli's theorem, $\overline{W}_{\Sigma}$ is a compact Lie group. 
Let $\E(\Sigma)^{\overline{W}_{\Sigma}}$ be the space of smooth 
$\overline{W}_{\Sigma}$-invariant functions on~$\Sigma$.  
It follows from  Schwarz's theorem 
\cite{Schwarz} that there exist $f_1,\ldots,f_n\in \E 
(\Sigma)^{\overline{W}_{\Sigma}}$ such that
\[F^{*}\E(\mathbf{R}^{n})=\E (\Sigma)^{\overline{W}_{\Sigma}}, \]
where $F=(f_1,\ldots,f_n)$.
Finally, we use Theorem~\ref{teo-basic-forms} to extend each 
function $f_i$ to a basic function on $M$ denoted by the 
same letter and this finishes the proof.
\end{proof}

%%%%%%%%%%%%%%%%%%%%%%%%%%%%%%%%%%%%%%%%%%%%%%%%%%%%%%%%%%%%%%%%%%%%%%%%%%%%%%%%%%%%%%%%%%%%%%%%%%%%%%
%%%%%%%%%% SECAO APLICACAO TRANSNORMAL
%%%%%%%%%%%%%%%%%%%%%%%%%%%%%%%%%%%%%%%%%%%%%%%%%%%%%%%%%%%%%%%%%%%%%%%%%%%%%%%%%%%%%%%%%%%%%%%%%%%%%%%

\section{Proof of Theorem \ref{teo-s.r.f.s-transnormal}}

Throughout this section we assume that $\F$ is a s.r.f.s.\ on a 
complete simply connected riemannian manifold $M$ 
with compact leaves and a flat section $\Sigma$, and
we prove Theorem~\ref{teo-s.r.f.s-transnormal}.

%We start by recalling 

%%%%
%%%LEMA
%%%
\begin{lema}
\label{lemma-reflection}
Let $\Omega_0$ and $\Omega_1$ be connected components of 
$M_{r}\cap\Sigma$ such that $\partial\Omega_0\cap\partial\Omega_1$ 
contains a wall $N$. Let $\sigma_0$ be a local section that intersects 
$\Omega_0$. Then there exists a local section $\sigma_1$ 
which intersects $\Omega_1$ and an isometry 
$\varphi:\sigma_0\cup \Omega_0\rightarrow \sigma_1\cup\Omega_1$ such that $\varphi(x)\in L_{x}$ for all  
 $x\in \Omega_0\cup \sigma_0$. 
In particular, $\varphi$ coincides with each holonomy with source in 
$\sigma_0\cup \Omega_0$ and target in $\sigma_1\cup \Omega_1$.
\end{lema}
\begin{proof}
Let $\beta$ be a curve contained in a regular leaf such that $\beta(0)\in\Omega_0$ and $\beta(1)\in\Omega_1$. To begin with,
we want to extend the singular holonomy $\varphi_{[\beta]}$ to an isometry $\varphi:\Omega_0\rightarrow\Omega_1.$ 
Since $\Omega_0$ is convex, flat and simply connected 
(see Theorem \ref{teo-fundamental-domain}), 
there exists a unique vector $\xi$ such that $\exp_{\beta(0)}(\xi)=x$. 
Let $\xi(\cdot)$ be the normal parallel transport of the vector $\xi$ 
along the curve $\beta$. We define $\varphi(x)=\exp_{\beta(1)}(\xi(1))$. 
It follows from Theorem~\ref{frss-eh-equifocal} that 
$\exp_{\beta(t)}(\xi(t))\in L_{x}$, and hence $\varphi(x)\in L_x$.  
Since $\varphi(x)\in L_x$ and $\Omega_1$ is the interior of a 
fundamental domain, the map $\varphi$ restricted to a neighborhood 
of each regular point of $\Omega_0$ coincides with a regular holonomy. 
This implies that $\varphi$ is an isometry, and hence  
$\varphi$ is smooth.   

Next, we want to extend $\varphi$ so that it is
also defined on $\sigma_0$. Since the restriction of $\varphi$ to  
a neighborhood of each regular point of $\Omega_0$ coincides with a 
regular holonomy, it suffices to prove that $\sigma_0\cap\Omega_0$ 
has only one connected component. If $\sigma_0$ is centered at a 
regular point $q$, then the fact that it is contained in the 
slice of $q$ implies that it is a ball contained in 
$\Omega_0$. Suppose now that $\sigma_0$ is centered at a 
singular point $q.$ Since it is contained in the slice of $q$, 
it follows from Theorem \ref{sliceteorema} that the intersection 
of $\sigma_0$ with the singular stratification of $\Sigma_0$ 
has only one connected component and also that 
$q\in\partial\Omega_0.$ These facts together with the fact that 
$\overline{\Omega}_0$ is flat and simply connected imply that 
also in this case
$\sigma_0\cap\Omega_0$ has only one connected component. \end{proof}

%\claim The intersection of $\sigma$ with the singular stratification of $\Sigma$ has only one connected component and also that $q\in\partial\Omega_0.$

%Let $\Pi:\mathbf{R}^n\rightarrow\Sigma$ be the  riemannian  universal covering map of the section $\Sigma.$ Since $\sigma$ is simply connected we can lift it to an open set $\widetilde{\sigma}.$ Since $\overline{\Omega}_0$ is simply connected we can lift it to a set $\overline{\widetilde{\Omega}}_0$ so that $\overline{\widetilde{\Omega}}_0\cap \widetilde{\sigma}$ is not empty.  The above claim implies that the intersection of $\widetilde{\sigma}$ with the singular stratification of $\Sigma$ has only one connected component and also that $q\in\partial\Omega_0.$   

%This fact and the fact that $\sigma$ is simply connected implies that $\tilde{\sigma}$ is a ball centered at $\tilde{q}$

%Using the fact that $\sigma$ is contained in a slice, Theorem \ref{teo-fundamental-domain} and the fact that $\Sigma$ is flat, we can conclude that %$\sigma\cap\Omega_0$ has only one connected component. Therefore we can extend $\varphi$ to be also defined on $\sigma$. 

%%%%%%%%%%%%%%%%%%%%%%
%%%%%%%%%%%%%%%%%%% NOVA PROP RECOBRIMENTO EH SUBGROUP DE COXETER
%%%%%%%%%%%%%%%%%%
\begin{prop}\label{cov-transf}
Let $\Pi:\mathbf{R}^{n}\rightarrow\Sigma$ be the riemannian 
universal covering map of the 
section $\Sigma$. Then: 
\begin{enumerate}
\item[(a)] There exists a locally finite family $\planes$ of hyperplanes in $\mathbf{R}^{n}$ which is invariant under the action of the group of isometries $\widetilde{W}$ generated by the reflections in the hyperplanes of $\planes$. In addition, the projection $\Pi (\planes )$ is the singular stratification on $\Sigma$.
\item[(b)] The group of
covering transformations of $\Pi$ is a subgroup of $\widetilde{W}.$ 
\item[(c)] Each singular holonomy $\varphi_{[\beta]}\in W_{\Sigma}$  is a 
restriction of a local isometry given by the 
composition of a finite number of reflections in the walls of the singular 
stratification of $\Sigma.$  
\end{enumerate}
\end{prop}
\begin{proof}
(a) The existence of the locally finite family $\planes$ 
of hyperplanes in  
$\mathbf{R}^{n}$ follows from the facts that 
$\Pi$ is a covering map and the singular stratification 
in the section $\Sigma$ is locally finite. 
We still need to prove that $\planes$ is invariant under the action of 
the group of isometries $\widetilde{W}$.
Let $\widetilde{H}_{0}$ be an hyperplane in $\planes$, and let 
$\widetilde{w}$ be the reflection in $\widetilde{H}_{0}$. 
Given an hyperplane $\widetilde{H}_1$, we want to prove that 
$\widetilde{w}(\widetilde{H_1})$ is an hyperplane in $\planes.$      

Let $\widetilde{\gamma}$ be the segment of line that joins 
a point $\widetilde{\gamma}(0)\in\widetilde{H}_0$ to a point 
$\widetilde{\gamma}(1)\in \widetilde{H}_1$ such that 
$\widetilde{\gamma}$ is orthogonal to $\widetilde{H}_1$ at 
$\widetilde{\gamma}(1)$.
Let us define $\gamma$ as the geodesic segment $\Pi(\widetilde{\gamma}).$ 
Then we cover $\widetilde{\gamma}$ (respectively, $\gamma$) 
by neighborhoods $\widetilde{U}_0,\ldots,\widetilde{U}_n$ 
(respectively, by neighborhoods 
 $U_0,\ldots,U_n$) so that $\Pi:\widetilde{U}_i\rightarrow U_i$ 
is an isometry, $\widetilde{\gamma}(0)\in \widetilde{U}_0$ and 
$\widetilde{\gamma}(1)\in \widetilde{U}_n$. 
 
Define the singular holonomy $\varphi_0: U_0\rightarrow U_0$ 
such that 
$\varphi_0\Pi|_{\widetilde{U}_{0}}=\Pi \widetilde{w}|_{\widetilde{U}_{0}}$.
By induction, we define a singular holonomy 
$\varphi_{n}:U_{n}\rightarrow U_{-n}$ such that 
$\varphi_{n-1}|_{U_{n-1}\cap U_{n}}=\varphi_{n}|_{U_{n-1}\cap U_{n}}$. 
Owing to 
$\varphi_{n-1}\Pi|_{\widetilde{U}_{n-1}} =\Pi \widetilde{w}|_{ \widetilde{U}_{n-1}}$, 
we conclude that
\begin{eqnarray}
\label{eq-1-prop-cov-transf}
\varphi_{n}\Pi|_{\widetilde{U}_{n}} =\Pi \widetilde{w}|_{ \widetilde{U}_{n}}.
\end{eqnarray}
Let $H_1$ be a wall whose closure contains $\gamma(1)$ and such that 
$H_1\subset \Pi(\widetilde{H}_1)$. Since $\varphi_n$ is a singular holonomy, $\varphi_n(H_1)$ is contained in the singular stratification of $\Sigma$. 
This fact together with Equation~(\ref{eq-1-prop-cov-transf}) yield that
$\widetilde{w}(\widetilde{H}_1)\in \planes$.

(b) Let $\gamma$ be a loop so that $\gamma(0)=x_0=\gamma(1)$, and 
$\tilde{\gamma}$ be the lift of $\gamma$ such that 
$\tilde{\gamma}(0)=\tilde{x}_0$.  Without loss of generality, 
we may assume that $\gamma$ meets the 
singular stratification only in the walls and always tranversally to them.

Let $0=t_0<\cdots < t_{n+1} =1$ be a partition such that
\begin{enumerate}
\item[(i)] $\gamma|_{(t_{i-1},t_{i})}$ has only regular points,
\item[(ii)] $\gamma(t_i)$ belongs to a wall for $0<i<n+1$.
\end{enumerate}
By induction define $x_i$ (respectively $\tilde{x}_i$) as the reflection 
of $x_{i-1}$ (respectively $\tilde{x}_{i-1}$) in the wall that contains 
$\gamma(t_i)$ (respectively $\tilde{\gamma}(t_i)$). 
Lemma \ref{lemma-reflection} implies that $ x_i\in L_{x_0}$ and 
$\Pi(\tilde{x}_i)=x_i.$
By construction, $\tilde{x}_n$ and $\tilde{\gamma}(1)$ both belong to the 
same Weyl chamber $\tilde{\Omega}.$ Note that 
$\Pi:\tilde{\Omega}\rightarrow\Omega$ is a diffeomorphism, where $\Omega$ is the connected 
component of $M_{r}\cap\Sigma$ that contains $x_0.$ Since $\overline{\Omega}$ 
is a fundamental domain, $L_{x_0}$ meets $\Omega$ only at $x_0.$ We 
conclude that $\Pi(\tilde{x}_n)=x_0=\Pi(\tilde{\gamma}(1)) .$
Since $\Pi:\tilde{\Omega}\rightarrow\Omega$ is a diffeomorphism, we 
have that $\tilde{x}_n=\tilde{\gamma}(1)$.
Therefore $\tilde{\gamma}(1)=\tilde{x}_n=g_n\cdots g_1 \cdot 
\tilde{x}_0=g_n\cdots g_1 \cdot\tilde{\gamma}(0)$, where $g_i$ is a reflection in the wall that 
contains $\tilde{\gamma}(t_i)$. In other words, we conclude that the covering 
transformation that sends $\tilde{\gamma}(0)$ to $\tilde{\gamma}(1)$ is 
$g_n\cdots g_1$. 

(c) Let $\gamma$ be a curve in $\Sigma$ such that $\gamma(0)=\beta(0)$ and $\gamma(1)=\beta(1).$ As in item~(b), we choose a 
partition $0=t_0<\cdots < t_{n+1} =1$  such that
$\gamma|_{(t_{i-1},t_{i})}$ has only regular points, and
$\gamma(t_i)$ belongs to a wall for $0<i<n+1$.

%$0<t_1<\cdots < t_n<1$  so that
% $\gamma|_{(t_i,t_{i+1})}$ has only regular points and $\gamma(t_i)$ belongs to a wall. 
 
Let $w_i$ be the singular holonomy that is the reflection in the wall that contains $\gamma(t_i)$. According to 
Lemma \ref{lemma-reflection} we can define the source of each $w_i$ so that 
$w_n\cdots w_1: U\rightarrow \Sigma$ is well defined, where $U$ is the source of $\varphi_{[\beta]}.$
% \begin{enumerate}
% \item[1)] the source of $w_i$ contains the segment $\gamma|_{[t_{i-1},t_i]},$
% \item[2)] $w_n\cdots w_1: U\rightarrow \Sigma$ is well defined, where $U$ is the source of $\varphi_{[\beta]}.$
% \end{enumerate}
We want to prove that 
  \begin{eqnarray}
  \label{Eq-prop-cov-transf-c}
 \varphi_{[\beta]}=w_n\cdots w_1|_{U}.
 \end{eqnarray}
For $i=0$ and $i=1$, define $\Omega_i$ to be 
the connected component of $M_r\cap\Sigma$ that contains $\beta(i)$. 
Since $\Omega_1$ is the interior of a fundamental domain, 
$L_{\beta(0)}\cap \Omega_1=\{\beta(1)\}$. This fact together with the 
properties of singular holonomies imply that 
$w_n\cdots w_1\beta(0)=\beta(1)$. 
We conclude that $(w_n\cdots w_1)^{-1}\varphi_{[\beta]}$ is an 
holonomy that fixes $\beta(0)$. 
Since $\Omega_{0}$ is the interior of a 
fundamental domain, we get that 
$(w_n\cdots w_1)^{-1}\varphi_{[\beta]}=I$, 
where $I$ is the identity with germ at $\beta(0)$, and this 
implies Equation~(\ref{Eq-prop-cov-transf-c}).
\end{proof}

%%%%%%%%%%%%%%%%%%%%%%%%%%%%%%%%%%%%%%%%%%%%%%%%%%%%%%%%%%%%%%%%%%%%%%%%%%%%%%%%%%%%%%%%%%%%%%%%%%%

%%%%%%%%%%%%%%%%%%%%%%%%%%%%%%%%%%%%%%%%%%%%%%%%%%%%%%%%%%%%%%%%%%%%%%%%%%%%%%%%%%%%%%%%%%%%%%%%%%%%%%%%%%%%%%%%%%%
%%%PROP DE TERNG E BOURBAKI
%%%%%%%%%%%%%%%%%%%%%%%%%%%%%%%%%%%%%%%%%%%%%%%%%%%%%%%%%%%%%%%%%%%%%%%%%%%%%%%%%%%%%%%%%%%%%%%%%%%%%%%%%%%%%%%%%%%%

It follows from Terng~\cite[App.]{Te1}~that~$\widetilde{W}$ 
is a \emph{Coxeter group}, i.e.~the 
subgroup of isometries $\widetilde{W}$ is generated by reflections,  
the topology induced in $\widetilde{W}$ from the group of isometries of 
$\mathbf{R}^n$ is discrete and the action on $\mathbf{R}^n$ is proper. 
Since $\planes$ is invariant by the action of $\widetilde{W}$, we have the 
following results (see Bourbaki~\cite{Bou}, Ch.~V \S3 Propositions~6,
7, 8, and~10, and Remarque~1 on~p.86; Ch.~VI \S2 Proposition~8 and 
Remarque~1 on~p.180). 

\begin{prop}
\label{prop-Bourbaki}\label{bourbaki}
Let $\widetilde{W}$ be the Coxeter group defined in 
Proposition~\ref{cov-transf}. Then:
\begin{enumerate}
\item We have that $\widetilde{W}$ is a direct product of irreducible Coxeter groups $\widetilde{W}_i$ $(1 \leq i \leq s )$ and,
after a possible adjustment of the origin, there exists a decomposition of $\mathbf{R}^n$ into an
orthogonal direct sum of subspaces $E_i$ $(0 \leq i \leq s)$ 
such that $\widetilde{w} (x_0,x_1,\ldots,x_s)=(x_0,\tilde{w}_1(x_1),\ldots,\tilde{w}_s(x_s))$ 
for $\tilde{w}=\tilde{w}_1\cdots \tilde{w}_s \in \widetilde{W}$ with $\tilde{w}_i\in \widetilde{W}_i.$ 
\item Let $\planes_i$ be the set of hyperplanes of $E_i$ whose  reflections generate $\widetilde{W}_i.$ Then the set $\planes$ consists of hyperplanes of the form
\[H=E_0\times E_1\times\cdots\times E_{i-1}\times H_i \times E_{i+1}\times\cdots \times E_s\]
with $H_i \in \planes_i$ and $i=1,\ldots,s$.
\item Every chamber $C$ is of the form $E_0\times C_1\times \ldots\times C_s,$ where, for each~$i$, the set $C_i$ is a chamber defined in $E_i$ by the set of hyperplanes $H_i.$
\item Each $C_i$ is an open simplicial cone or an open simplex if $\widetilde{W}_i$ is finite or infinite respectively.
\item If $\widetilde{W}_i$ is infinite, then it is an 
\emph{affine Weyl group},
i.e.~there exists a unique root system $\Delta_i$ in $E_i$ such that 
$\widetilde W_i$ can be written as the semidirect product 
$\overline{W_i}\ltimes\Gamma_i$, where 
$\overline{W_i}$ is the Weyl group associated to $\Delta_i$ 
and $\Gamma_i$ is a 
group of traslations of $E_i$ of rank equal to $\dim E_i$; 
\end{enumerate}
\end{prop}

%%%%%%%%%%%%%%%%%%%%%%%%%%%%%%%%%%%%%%%%%%%%%%%%%%%%%%%%%%%%%%%%%%%%%%%%%%%%%%%%%%%%%%%%%%%%%%%%%%%%%%%%%%%%%%%
%%%%%%%% PROP O GRUPO DE WEYL  ADMITE BASE DE HILBERT 
%%%%%%%%%%%%%%%%%%%%%%%%%%%%%%%%%%%%%%%%%%%%%%%%%%%%%%%%%%%%%%%%%%%%%%%%%%%%%%%%%%%%%%%%%%%%%%%%%%%%%%%%%%%%%%%%%
\begin{prop}\label{free-polyn-dense-subalg}
Let $i$ be an index such that
$\widetilde{W}_i$ is an affine Weyl group acting on $E_i$
as in item~(e) of Proposition~\ref{bourbaki}.  
Then $\mathcal E(E_i)^{\widetilde{W}_i}$ contains a free 
polynomial subalgebra $\mathcal A$ on $\dim E_i$ generators which is 
dense in the sup-norm.
\end{prop}
\begin{proof}

The following argument is extracted from~\cite{HOL},
Theorem~7.6 and Corollary~7.7. Since the index $i$ is fixed,
throughout the proof we will drop it from the notation; we also 
identify $E_i$ with $\R^m$.  
Now the affine Weyl group $\widetilde W$ is the semi-direct
product $\overline W\ltimes\Gamma$, where $\overline W$ is the isotropy 
subgroup at zero and $\Gamma$ is a lattice of translations
of $\R^m$. 
Here zero is assumed to belong to one hyperplane
from each family of parallel singular hyperplanes. 
Since $\Gamma$ is a normal subgroup of $\widetilde W$,
the algebra of $\widetilde W$-invariant smooth functions 
on $\R^m$ can be written
\[ \mathcal E(\R^m)^{\widetilde W}\cong (\mathcal 
E(\R^m)^\Gamma)^{\overline W}
\cong \mathcal E(\R^m/\Gamma)^{\overline W}. \]
Note that $\R^m/\Gamma$ is a compact torus.
Since $\overline W$ is a 
finite group, 
$\mathcal E(\R^m/\Gamma)^{\overline W}$ separates the $\overline 
W$-orbits 
in $\R^m/\Gamma$, and this implies that 
$\mathcal E(\R^m)^{\widetilde W}$ separates the
$\widetilde W$-orbits in $\R^m$. 

Denote by $\Gamma^*$ the dual lattice of $\Gamma$ 
in $\R^{m*}$.
For each $\gamma\in\Gamma^*$, 
$x\in\R^m\mapsto e^{2\pi\sqrt{-1}\gamma(x)}\in\C$
induces a complex-valued smooth 
function on $\R^m/\Gamma$ which is denoted by
$e^{2\pi\sqrt{-1}\gamma}$.
Let $\C[\Gamma]$ denote the complex algebra consisting 
of finite linear combinations 
with complex coefficients of elements of 
$\{e^{2\pi\sqrt{-1}\gamma}\}_{\gamma\in\Gamma^*}$,
and let $\R[\Gamma]$ denote the subalgebra 
of $\C[\Gamma]$ consisting of real-valued 
functions. 
Since $\{e^{2\pi\sqrt{-1}\gamma}\}_{\gamma\in\Gamma^*}$ is an 
orthogonal 
basis of $L^2(\R^m/\Gamma)$ consisting of
eigenfunctions of the Laplacian of $\R^m/\Gamma$, it
is known that 
$\C[\Gamma]$ is a dense subalgebra of $\mathcal E(\R^m/\Gamma)\otimes\C$
with respect to the sup-norm. By taking real parts 
and averaging, we get that 
$\mathcal A:=\R[\Gamma]^{\overline W}$ is a dense subalgebra of 
$\mathcal E(\R^m/\Gamma)^{\overline W}$. It remains
to prove that $\mathcal A$ is a 
free polynomial algebra on $m$ generators. 

The $\overline W$-action on $\mathcal E(\R^m/\Gamma)$
extends $\C$-linearly to $\mathcal E(\R^m/\Gamma)\otimes\C$,
and, it follows from~\cite{Bou}, 
Ch.~VI \S3 Th\'eor\`em~1, that $\C[\Gamma]^{\overline W}$ is a free 
polynomial algebra on $m$ generators. Moreover, as explained 
in that book, the generators can be chosen 
in a special way, as follows. 

Let $\Delta$ be the root system in $\R^{m*}$
associated to $\overline W$.
For each root $\alpha\in\Delta$, the corresponding
inverse root $\breve\alpha\in\R^m$ is defined 
as being  
$\breve\alpha=2h_\alpha/||h_\alpha||^2$, where 
$h_\alpha$ is the element of $\R^m$ satisfying 
$\inn{h_\alpha}{x}=\alpha(x)$ for all $x\in\R^m$.
Identifying the translations of $\Gamma$ with 
elements of $\R^m$, we have that the lattice of 
inverse roots 
coincides with $\Gamma$~\cite[Ch.~VI \S2 Proposition~1]{Bou}.
%~\cite[Proposition~7.9]{BD}. 
It follows that $\Gamma^*$ coincides
with the dual lattice of the lattice of 
inverse roots, which is
by definition the lattice of weights. Chosen a Weyl chamber,
there is a distinguished basis $\gamma_1,\ldots,\gamma_n$
of $\Gamma^*$ whose elements are called the fundamental weights. 
The $\overline W$-action on $\R^{m*}$ permutes the elements
of $\Gamma^*$, and, plainly, 
$w\cdot e^{2\pi\sqrt{-1}\gamma} = e^{2\pi\sqrt{-1}(w\cdot\gamma)}$
for $w\in\overline W$.
The averaging operator
is the $\C$-linear map $S:\C[\Gamma]\to\C[\Gamma]^{\overline W}$
given by 
$S(e^{2\pi\sqrt{-1}\gamma})=\frac1{|\overline W|}\sum_{w\in \overline 
W}e^{2\pi\sqrt{-1}(w\cdot\gamma)}$.
Then free generators $x_1,\ldots,x_n$ of $\C[\Gamma]^{\overline W}$ 
can be taken to be $x_i=S(e^{2\pi\sqrt{-1}\gamma_i})$. 

Finally, we construct free generators for 
$\mathcal A$. Note that 
$-\gamma_i$ and $\gamma_i$ belong to the same
$\overline W$-orbit if and only if 
$-\gamma$ and $\gamma$ belong to the same $\overline W$-orbit for all 
$\gamma$ in the $\overline W$-orbit of $\gamma_i$, and this 
holds if and only if $x_i$ is real-valued. In general, 
according to~\cite{BD}, Note~4.2, there is an involution $\varrho$
of the set $\{1,\ldots,n\}$ such that 
$-\gamma_i\in \overline W(\gamma_{\varrho(i)})$. By rearranging
the indices, we may thus assume that $\varrho(i)=i$ for
$i=1,\ldots, p$ and $\varrho(p+i)=p+q+i$ for $i=1,\ldots,q$, where
$p+2q=n$. It follows that $x_1,\ldots,x_p$ are real-valued,
and $x_{p+i}$, $x_{p+q+i}$ are complex-conjugate
for $i=1,\ldots,q$. Therefore real-valued free generators 
$y_1,\ldots,y_n$ of $\mathbf C[\Gamma]^{\overline W}$ can be chosen
so that $y_i=x_i$ for $i=1,\ldots,p$ and 
$y_{p+i}=\Re x_{p+i}$, $y_{p+q+i}=\Im x_{p+q+i}$ for 
$i=1,\ldots,q$. It is clear that $y_1,\ldots,y_n$ generate 
$\mathcal A$ as an algebra over $\R$. 
\end{proof}

%%%%%%%%%%%%%%%%%%%%%%%%%%%%%%%%%%%%%%%%%%%%%%%%%%%%%%%%%%%%%%%%%%%%%%%%%%%%%%%%%%%%%%%%%%%%%%%%%%%%%%%%%%%%%%%%%%%%%%
%%%%%%%%%%%%%%%%%%
%%%%%%%%%%%%%%%%% PROP EXISTE UMA APLICACAO W INVARIANTE NO ESP EUCLIDIANO
%%%%%%%%%%%%%%%%%%%%%%%%%%%%%%%%%%%%%%%%%%%%%%%%%%%%%%%%%%%%%%%%%%%%%%%%%%%%%%%%%%%%%%%%%%%%%%%%%%%%%%

\begin{prop}
\label{prop-H-invariante}
Let $\widetilde{W}$ be the Coxeter group defined in Proposition \ref{cov-transf}. Then there exists a map $\widetilde{F}:\mathbf{R}^n\rightarrow\mathbf{R}^n$ which is $\widetilde{W}$-invariant and separates the $\widetilde{W}$-orbits.
\end{prop}
\begin{proof}
As remarked in Proposition \ref{prop-Bourbaki}, the Coxeter group 
$\widetilde{W}$ is the direct product of irreducible Coxeter groups 
$\widetilde{W}_i$ with $1\leq i\leq s$. 
If $\widetilde{W}_i$ is finite,  follows from Chevalley~\cite{Chevalley} 
that there exists a map $\widetilde{F}_i:E_i\rightarrow E_i$ which is 
$\widetilde{W}_i$ invariant and separates the $\widetilde{W}_i$-orbits. 
On the other hand, if $\widetilde{W}_i$ is infinite,  
Proposition~\ref{free-polyn-dense-subalg} implies the existence of a 
map $\widetilde{F}_i$ with the same properties as above. Finally we 
define 
$\widetilde{F}: E_0\oplus E_1\oplus\cdots\oplus E_s\rightarrow E_0\oplus E_1\oplus\cdots\oplus E_s$ to be
$F(x_0,x_1,\ldots,x_s)=(x_0,\widetilde{F}_{1}(x_1),\ldots,\widetilde{F}_{s}(x_s))$.
\end{proof}

%%%%%%%%%%%%%%%%%%%%%%%%%%%%%%%%%%%%%%%%%%%%%%%%%%%%%%%%%%%%%%%%%%%%%%%%%%%%%%%%%%
%%%%%%%%%%%%%%%%%%% PROVA DO TEOREMA APLICACOES TRANSNORMAIS
%%%%%%%%%%%%%%%%%%%%%%%%%%%%%%%%%%%%%%%%%%%%%%%%%%%%%%%%%%%%%%%%%%%%%%%%%

\medskip

\noindent\textit{Proof of Theorem~\ref{teo-s.r.f.s-transnormal}.}
It follows from Proposition~\ref{cov-transf}(b) and 
Proposition~\ref{prop-H-invariante} that there exists a map 
$\hat{F}:\Sigma\rightarrow\mathbf{R}^{n}$ such that 
$\hat{F}\circ\Pi=\widetilde{F}$. 
Now Proposition~\ref{prop-H-invariante} and Proposition~\ref{cov-transf}(c) 
imply that the 
map~$\hat{F}$ is~$W_{\Sigma}$-invariant, 
where $W_{\Sigma}$ is the Weyl pseudogroup of $\Sigma$. 
We can apply Theorem~\ref{teo-basic-forms} to extend the map 
$\hat{F}$ to a map $F: M\rightarrow\mathbf{R}^{n}$ such that 
the leaves of $\F$ 
coincide with the level sets of $F$.  
Finally, the transnormality 
of $F$ follows from the fact that the regular leaves of 
$\F$ form a riemannian foliation with sections
(see Molino~\cite[p.~77]{Molino}), and this 
finishes the proof of the theorem.\eop

%%%%%%%%REFERENCIAS
%%%%%%%%%%%%%%
\bibliographystyle{amsplain}

\end{document}